\input amstex 
\documentstyle{amsppt}
\input bull-ppt
\keyedby{bull380/lic}
\define\ocl{\operatorname{cl}}
\define\Aut{\operatorname{Aut}}
\define\gl{\operatorname{GL}}

\topmatter
\cvol{28}
\cvolyear{1993}
\cmonth{April}
\cyear{1993}
\cvolno{2}
\cpgs{315-323}
\title Zariski Geometries \endtitle
\author Ehud Hrushovski and Boris Zilber\endauthor
\shortauthor{Ehud Hrushovski and Boris Zilber}
\shorttitle{Zariski Geometries}
\address Massachusetts Institute of Technology and Hebrew
University at Jerusalem\endaddress
\cu Department of Mathematics, Massachusetts Institute of 
Technology,
Cambridge, Massachusetts 02139\endcu
\ml ehud\@math.mit.edu\endml
\address Kemerovo University and University of Illinois at 
Chicago\endaddress
\cu Department of Mathematics, Kemerovo University, 
Kemerovo,
650043, Russia\endcu
\ml zilber\@kemucnit.kemerovo.su\endml
\date February 12, 1992 and, in revised form, August 24, 
1992\enddate
\thanks The first author was supported by NSF grants DMS 
9106711 and
DMS 8958511\endthanks
\abstract We characterize the Zariski topologies over an 
algebraically
closed field in terms of general dimension-theoretic 
properties.
Some applications are given to complex manifold and to 
strongly minimal
sets.\endabstract
\endtopmatter

\document

\heading 1. Introduction\endheading
There is a class of theorems that characterize certain 
structures by their
basic topological properties. For instance, the only 
locally compact connected
fields are $\bold R$ and $\bold C$. These theorems refer 
to the classical
topology on these fields. The purpose of this paper is to 
describe a
similar result phrased in terms of the Zariski topology. 

The results we offer differs from the one considered above 
in that we
do not assume in advance that our structure is a field or 
that it carries an
algebraic structure of any kind. The identification of the 
field structure is
rather a part of the conclusion.

Because the Zariski topologies on two varieties do not 
determine the Zariski
topology on their product (and indeed the topology on a 
one-dimensional variety
carries no information whatsoever), the data we require 
consists not only
of a topology on a set $X$, but also of a collection of 
compatible topologies
on $X^n$ for each $n$. Such an object will be referred to 
here as a geometry.
It will be called a {\it Zariski geometry\/} if a 
dimension can be assigned to
the closed sets, satisfying certain conditions described 
below. Any smooth
algebraic variety is then a Zariski geometry, as is any 
compact complex
manifold if the closed subsets of $X^n$ are taken to be 
the closed holomorphic
subvarieties. 

If $X$ arises from an algebraic curve, there always exist 
large
families of closed subsets of $X^2$; specifically, there 
exists a family of
curves on $X^2$ such that through any two points there is 
a curve in the family
passing through both and another separating the two. An 
abstract Zariski
geometry with this property is called {\it very ample}. By 
contrast, there are
examples of analytic manifolds $X$ such that $X^n$ has 
very few closed analytic
submanifolds and is not ample. Precise definitions will be 
given in the next
section; the complex analytic case is discussed in \S4. 
Our main result is 

\thm{Theorem 1} Let $X$ be a very ample Zariski geometry. 
Then there exists a
smooth curve $C$ over an algebraically closed field $F$, 
such that $X$, $C$ are
isomorphic as Zariski geometries. $F$ and $C$ are unique, 
up to a field
isomorphism and an isomorphism of curves over $F$.\ethm

A weaker version is available in higher dimensions (and 
actually follows easily
from Theorem 1). If $X$ is a higher-dimensional Zariski 
geometry and there
exists a family of curves on $X$ passing through any two 
points and separating
any two, then there exists a dense open subset of $X$ 
isomorphic to a 
variety. Globally, $X$ probably arises from an algebraic 
space in the sense of
\cite A, but we do not prove this. In this survey we will 
concentrate on the 
one-dimensional case (though, of course, we must consider 
the geometry of
arbitrary powers). 

\heading 2. Zariski Geometries\endheading
Zariski geometries will be defined as follows. Recall that a
topological space is {\it Noetherian\/} if it has the 
descending chain 
condition on closed subsets. (See \cite{Ha}.) A closed set 
is {\it
irreducible\/} if it is not the union of two proper closed 
subsets.
If $X$ is Noetherian, then every closed set can be written 
as a finite
union of irreducible closed sets. These are uniquely 
determined 
(provided that no one is a subset of the other) and are 
called the
{\it irreducible components\/} of the given set. We say 
that $X$ has dimension
$n$ if $n$ is the maximal length of a chain of closed 
irreducible sets 
$C_n\supset C_{n-1}\supset\cdots\supset C_0$. 

We will use the following notation: if $C\subseteq 
D_1\times D_2$,
$a\in D_1$, we let $C(a)=\{b\in D_2\:(a,b)\in C\}$.

\dfn{Definition 1} A Zariski geometry on a set $X$ is a 
topology on 
$X^n$ for each $n$ that satisfies the following:
\roster
\item"(Z0)" Let $f_i$ be a constant map $(f_i(x_1,\dots, 
x_n)=c)$ or a
projection $(f_i(x_1,\dots, x_n)=x_{j(i)})$. Let 
$f(x)=(f_1(x),\dots,
f_m(x))$. Then $f\: X^n 
\to X^m$ is continuous. The diagonals $x_i=x_j$ of $X^n$ 
are closed.
\item"(Z1)" Let $C$ be a closed subset of $X^n$, and let 
$\pi$ be the
projection to $X^k$. Then there exists a proper closed 
subset $F$ of
$\ocl(\pi C)$ such that $\pi C\supseteq\ocl(\pi C)-F$. 
\item"(Z2)" $X$ is irreducible and {\it uniformly\/} 
1-{\it dimensional\/}:
if $C\supseteq X^n \times X$ is closed, then for some $m$, 
for all
$a\in X^n$, $C(a)=X$ or $|C(a)|\le m$.
\item"(Z3)" (Dimension Theorem) Let $U$ be a closed 
irreducible subset
of $X^n$, and let $\Delta_{ij}$ be the diagonal $x_i=x_j$. 
Then every
component of $U\cap \Delta_{ij}$ has dimension $\ge\dim 
(U)-1$. 
\endroster

\enddfn

Axiom (Z0) includes the obvious  compatibility 
requirements on the
topologies.

$X$ is called {\it complete\/} (or {\it proper\/}) if all 
projection
maps are closed. Note that (Z1) is a weakening of 
completeness. We
prefer not to assume completeness axiomatically because we 
do not wish to 
exclude affine models. 

(Z2) is the assumption of one-dimensionality that we make 
in this paper.
The results still have consequences in higher dimensions, 
as will be 
seen in \S4.

Axiom (Z3) is valid in smooth algebraic (or analytic) 
varieties.
One still obtains information on other varieties by 
removing the
singular locus, applying the theorem, and going back. It 
is partly 
for this reason that it is important not to assume 
completeness
in (Z1).

We now turn to the ``ampleness'' conditions. By a {\it 
plane curve\/}
over $X$ we mean an irreducible one-dimensional subset of 
$X^2$. A
{\it family of plane curves\/} consists of a closed 
irreducible set
$E\subseteq X^n$ (parametrizing the family) and a closed 
irreducible
$C\subseteq E\times X^2$, such that $C(e)$ is a plane 
curve for generic
$e\in E$. 

\dfn{Definition 2} A Zariski geometry $X$ is {\it very 
ample\/} if there
exists a family $C\subseteq E\times X^2$ of plane curves 
such that:
\roster
\item"(i)" For generic $a,b\in X^2$ there exists a curve 
$C(e)$
passing through $a,\ b$. 
\item"(ii)" For any $a,b\in X^2$ there exists $e\in E$ 
such that
$C(e)$ passes through just one of $a,b$.
\endroster
If only (i) holds, then $X$ is called {\it ample.}
\enddfn

Axioms (Z0)--(Z3) alone allow a degenerate geometry, in 
which the 
closed irreducible subsets of $X^n$ are just those defined 
by equations of
the form $x_i=a_i$. More interesting are the {\it 
linear\/} geometries,
where $X$ is an arbitrary field (or division ring) and the 
closed subsets
of $X^n$ are given by linear equations. These are {\it 
nonample Zariski
geometries.} A thorough analysis of nonample Zariski 
geometries can
be carried out; they are degenerate or else closely 
related to the
above linear example. We will not describe this analysis 
here; see
\cite{HL}.

Theorem 1 describes the very ample Zariski geometries, 
while the 
nonample Zariski geometries were previously well 
understood. This leaves
a gap---the ample but not very ample Zariski geometries. 
We can show:

\thm{Theorem 2} Let $D$ be an ample Zariski geometry. Then 
there exists an
algebraically closed field $K$ and a surjective map $f\: 
D\to \bold P^1(K)$.
$f$ maps constructible sets to \RM(algebraically\/\RM) 
constructible sets\,\RM;
in fact off a certain finite set, $f$ induces a closed, 
continuous map on each
Cartesian power.\ethm

This represents $D$ as a certain branched cover of $\bold 
P^1$. Among
complex analytic manifolds, all such covers are algebraic 
curves. This is not
true in the Zariski context; there are indeed ample, not 
very ample Zariski
geometries, which do not arise from algebraic curves. We 
construct these as
formal covers of $\bold P^1$, starting from any nonsplit 
finite extension $G$
of a subgroup of $\Aut(\bold P^1)$ that cannot be realized 
as the automorphism
group of any algebraic curve. We obtain a finite cover of 
$\bold P^1$, with an
action of $G$ on it, and define the closed sets so as to 
include the pullbacks
of the Zariski closed subsets of $\bold P^1$ and the 
graphs of the
$G$-operators. The automorphism group of the cover is then 
too large to arise
from an algebraic curve.

The following is implicit in Theorem 1, but we wish to 
isolate it:

\thm{Theorem 3} Let $K$ be an algebraically closed field, 
and let $X$ be a
Zariski geometry on $\bold P^1 (K)$, refining the usual 
Zariski geometry.
Then the two geometries coincide.\ethm

These results sprang from two sources, which we proceed to 
describe. 

\heading 3. Strongly minimal sets\endheading
The original goal was to characterize an algebraically 
closed field $F$ in
terms of the collection of constructible subsets of $F^n$, 
rather than in terms
of the closed subsets. The motivation was the importance 
in model theory of
certain structures, called strongly minimal sets. (In 
particular they form the
backbone of any structure categorical in an uncountable 
power (see \cite{BL}).)

\dfn{Definition 3} A {\it structure\/} is an infinite set 
$D$ together with a
collection of subsets of $D^n\ (n=1,2,\dots)$ closed under 
intersections,
complements, projections and their inverses, and 
containing the diagonals. 
These are called the 0-{\it definable sets.} $D$ is {\it 
strongly minimal\/}
if it satisfies:
\roster
\item"(SM)" For every 0-definable $C\subseteq D^{m+1}$ 
there exists an integer
$n$ such that for all $a\in D^n$, letting $C(a)=\{b\in 
D\:(a,b)\in C\}$,
either $|C(a)|\le n$ or $|D-C(a)|\le n$. 
\endroster

\enddfn

These axioms can be viewed as analogs of (Z1), (Z2) for 
the class of
constructible sets. Under these assumptions one can prove 
the existence of a
well-behaved dimension theory. In particular, one can 
state an axiom analogous
to ampleness (non-local-modularity; strongly minimal sets 
{\it not\/}
satisfying it are well understood). See \cite{Z, HL}. 
However, it is not clear how to
state (Z3) in terms of constructible sets alone. In the 
absence of such an
axiom, it was shown in \cite{Hr1} that the analog of 
Theorem 2 is false and in 
\cite{Hr2} that the analog of Theorem 3 also fails. 

We note that Macintyre (see \cite{Mac}) characterized the 
strongly minimal
{\it fields\/} as the algebraically closed fields. (These 
are strongly minimal
sets with a definable field structure, i.e., the graphs of 
addition and
multiplication are 0-definable subsets of $D^3$.) A 
conjecture by Cherlin and
the second author that simple groups definable over 
strongly minimal sets are
algebraic groups over an algebraically closed field 
remains open.

\heading 4. Complex manifolds\endheading
Let $X$ be a (reduced, Hausdorff) compact complex analytic 
space, and consider
the topology $An$ on $X^n$ whose closed sets are the 
closed analytic 
subvarieties of $X^n$. Remmert's theorem then implies that 
the projection of a
closed set is closed. It can be shown further that if $U$ 
is a {\it locally
closed\/} subset of $X^{n+1}$, i.e., the difference of two 
$An$-closed
sets, then the projection of $U$ to $X^n$ is itself a 
finite union of locally
closed sets, or an $An$-{\it constructible set.} Thus the 
$An$-constructible
sets form a {\it structure\/} in the sense of \S3. If $V$ 
is a minimal analytic
subvariety of $X$, then $V$ with this structure is 
strongly minimal; further,
if one removes the singular locus of $V$, one obtains a 
Zariski geometry in the
sense of \S1. In this section we discuss some consequences 
of this observation.

In the analytic context, Theorem 2 resembles Riemann's 
existence theorem (the
part stating that a compact complex manifold of dimension 
1 is a finite cover
of the projective line). Indeed Riemann's existence 
theorem would follow from
Theorem 2, but the hypothesis of Theorem 2 includes 
ampleness, which we do not
know how to prove directly. Instead we offer a variation 
in dimensions $\ge 2$.
Note that the assumption on $M$ is true of a generic 
complex torus of dimension
$\ge 2$ (as will also follow from Proposition 3).

\thm{Proposition 1} Suppose $M$ is a compact K\"ahler 
manifold of complex
dimension $\ge 2$, with no proper infinite analytic 
subvarieties. If
$H_1(M)\ne 0$, then $M$ is a complex torus.\ethm

This follows from Theorem 2, as follows. If $M$ were 
ample, then it would be
a finite branched cover of $\bold P^1(K)$ for some field 
$K$, which in the
present context also has the structure of a complex 
analytic space. $K$ must
also have complex dimension $\ge 2$, which contradicts the 
classification of
the connected locally compact fields cited in the 
Introduction. Thus $M$ cannot
be ample. As mentioned in the previous section, this gives 
a strong ``Abelian''
condition; it is shown in \cite{HP} (in a much more 
general context) that if
ampleness fails in a given geometry, then every closed 
irreducible subset of a
group $A$ supported by that geometry must be a coset of a 
closed subgroup.
We apply this to the Albanese variety $A$ of $M$ 
(\cite{GH, p.\ 331}). A fiber
of the map from $M$ to $A$ must be finite, or we violate 
the assumption on
infinite analytic subvarieties. The image of $M$ in $A$ is 
closed and
irreducible, hence a coset of a closed subgroup $S$ of 
$A$; ``closed'' here
means a complex analytic subvariety, hence a subtorus, of 
$A$. By translation
we obtain a finite holomorphic map $f\: M\to S$. The 
branch locus of this map
gives an analytic subvariety locally of dimension 
$\dim(M)-1$ and hence must be
empty. Thus $M$ is a finite covering of the complex torus 
$S$ and hence is
itself a complex torus.

Theorem 3 resembles Chow's theorem that a closed analytic 
subvariety of
projective space must be algebraic (see \cite{GH}). Indeed 
a somewhat more
general result can easily be deduced from the theorem.

\thm{Proposition 2} Let $X$ be a complex algebraic 
variety. View $X$ as a
complex analytic space \RM(as in \cite{HA, \RM{Appendix 
1}}\RM). Then any closed
analytic subvariety of $X$ is an algebraic subvariety.\ethm

This is easily seen by working with constructible sets, 
defined to be finite
Boolean combinations of closed sets. We have two Zariski 
geometries on $X$,
given by the algebraic and the analytic structures. $X$ 
contains smooth
algebraic curves, on which again two geometries are 
induced. By Theorem 1, one
obtains the same geometry on each such curve. It follows 
easily that the two
geometries have the same class of constructible sets, so 
any closed analytic
subvariety $V$ of $X$ is algebraically constructible. Then 
it is easy to see
that $V$ must in fact be Zariski closed.

Another application of \cite{HP} yields the following 
statement:

\thm{Proposition 3} Let $G=C^n/L$ be a complex torus. Then 
either $G$ has closed
analytic subgroups $G_1\subseteq G_2$, such that $G_2/G_1$ 
is isomorphic to a
nontrivial Abelian variety, or the only complex analytic 
submanifolds of $G$
are subtori and finite unions of their cosets.\ethm

\heading 5. Method of proof\endheading
We begin with a description of the proof of Theorem 2. 

\subheading{\num{1.} A universal domain} We are given a 
Zariski geometry $X$,
which we think of as analogous to the Zariski geometry on 
a curve over an
algebraically closed field, and its powers. We wish to 
find an analog to the
finer notion of the $K$-topology (where $K$ is a subfield 
of an algebraically
closed field); in other words, we wish to have a concept 
of a closed set
``defined over a given substructure $K$''. This cannot be 
usefully done for
$X$ itself; we need to embed $X$ in a larger geometry 
$X^*$, a ``universal
domain''. (This is analogous to viewing a number field as 
a subfield of a larger
field of infinite transcendence degree.) This construction 
of a universal
domain is in fact a standard one in model theory; $X^*$ is 
obtained via the
compactness theorem of model theory or by using 
ultrapowers; see \cite{CK}
(saturated models) or \cite{FJ} (enlargements). 

\subheading{\num{2.} A combinatorial geometry} We assume 
(1) has been 
carried out and work directly with the universal domain 
$X^*$. An element $a$
of $X^*$ is {\it algebraically dependent\/} on a tuple of 
elements $b\in
X^{*n}$ if there exists a 0-definable closed set 
$C\subseteq X^{*n} \times
X^*$ such that $C(b)$ is finite and $a\in C(b)$. In the 
case of algebraically
closed fields, this coincides with the usual notion. It 
can be shown in general
that this notion of algebraic closure yields a {\it 
combinatorial pregeometry},
i.e., it satisfies the exchange axiom: If $a$ is 
algebraically dependent on
$b_1,\dots, b_n,c$ but not on $b_1,\dots, b_n$, then $c$ 
depends on $b_1,\dots,
b_n,a$. Formally, one can define a transcendence basis, 
dimension, etc.
In particular, the {\it rank\/} of a subset of $X^*$ is 
the size of any maximal
independent subset thereof.

\subheading{\num{3.} Group configurations} To construct a 
field, we will need
to find its additive and multiplicative groups. For this 
purpose we use a
general machinery (valid in the strongly minimal context 
and in fact
considerably beyond it) to recognize groups from the trace 
that they leave on
the combinatorial geometry. We will apply this machinery 
to the affine
translation group to obtain the field.

Suppose $G$ is a 1-dimensional group interpretable in the 
geometry $X^*$.
(Assume for simplicity that the elements of $G$ are points 
of $X^*$; the
graph of multiplication is assumed to be locally closed.) 
Let $a_1,a_2,\dots,
b_1,b_2,\dots$ be generic points of $X$, i.e., an 
independent set of elements
of $X^*$ (over some algebraically closed base substructure 
$B$). Let 
$c_{ij}=a_i+b_j$. Then $(c_{ij}\:i\in I,j\in J)$ forms an 
array of elements of
the combinatorial pregeom\-etry of $B$-dependence. It is 
easy to see that the
rank of any $m\times n$-rectangle in this array is $m+
n-1$. Moreover, for any
permutations $\sigma$ of $I$ and $\tau$ of $J$, the 
corresponding permutation
of the array arises from an automorphism of the geometry.

Conversely, suppose $(c_{ij}:i,j)$ is an array of elements 
of the combinatorial
pregeometry enjoying the above symmetry property and in 
which any $m\times
n$-rectangle has rank $m+n-1$. Then one proves the 
following theorem: There
exists a 1-dimensional Abelian group $G$ and independent 
generic elements
$a_i,b_j$ of $G$, such that $c_{ij}$ and $a_i+b_j$ depend 
on each other over
$B$. In particular, an infinite Abelian group is involved.

A similar theory is available for a connected group not 
necessarily Abelian, of
any dimension. We note in this connection Weil's theorem 
on ``group chunks''
\cite W. Weil's theorem allows the recognition of a group 
from generic data: A
binary function which is generically associative and 
invertible arises from a
definable group. Our result is similar but more general. 
In particular, it
permits the function $f$ to be multivalued (i.e., the 
graph of an algebraic
correspondence). (In other words, the hypothesis concerns 
algebraic rather than
rational dependence on certain generic points.) It is 
shown that given a trace
of associativity, there exists a group $G$ such that $f$ 
is conjugate (by a
multivalued correspondence between the given set and 
$G$\<) to the single-valued
function $xy^{-1}z$ of $G$. (See \cite{EH}.)

Applying the theorem to the two-dimensional group of 
affine transformations
of a field $F$, we obtain the following higher-dimensional 
analog:

\thm{Proposition 4} Let $c_{ij}\ (i,j=1,2,\dots)$ be a 
symmetric array of
elements of the dependence geometry over $B$. Suppose 
every $m\times n=$
rectangle of elements of $c_{ij}$ has rank $2m+n-2\ 
(m,n\ge 2)$. Then
there exists an algebraically closed field $F$ defined 
over $B$ and generic
independent elements $a_i,b_j,g_k$ of $F$, such that 
$c_{ij}$ and
$a_i+b_ig_j$ depend on each other.\ethm

We will merely use the existence of $F$. We note that if 
$c_{ij}=
a_i+b_i g_j$, then the elements $c_{ij}$ satisfy the 
relations: $(c_{ij}-
c_{ij'})/(c_{i'j}-c_{i'j'})=
(c_{ij}-c_{ij''})/(c_{i'j}-c_{i'j''})$. 

\subheading{\num{4.} Tangency} So far we have used only 
the strongly
minimal set structure, not the more detailed knowledge of
the identity
of the closed sets. This comes in via the notion of a 
specialization.
A map $f$ from a subset $C$ of $X^*$ to $X^*$  is said to 
be a
specialization (over $B$\<) if for every $B$-closed set 
$F\subset C^n$
and all $a_1,\dots, a_n\in C$, if $(a_1,\dots, a_n)\in F$ 
then $(fa_1,\dots,
fa_n)\in F$. Note that no such notion is available for $X$ 
itself. The
properness axiom immediately yields the extension theorem 
for 
specializations. The dimension theorem gives a more subtle 
property,
one of whose principal consequences is a ``preservation of 
number''
principle: If $C(a)$ is a closed set depending continuously 
on $a$,
$a\to a'$ is a specialization, and $C(a)$, $C(a')$ are both 
finite,
then there exists a specialization of $\{a\}\cup C(a)$ onto 
$\{a'\}
\cup C(a')$; in particular, the number of points in $C(a)$ 
cannot go
up (but remain finite) under specializations.

This allows us to formalize a notion of intersection 
multiplicity.
Let $C_1,C_2$ be curves in the ``surface'' $X\times X$, 
depending
as above on parameters $c_1,c_2$ (so $C_i=C^i(c_i)$ for some
$0$-closed
set $C^i$\<). Let $a_1,a_2$ be distinct points in the 
intersection of
$C_1,C_2$. If $(c_1,c_2,a_1,a_2)\to
(c'_1, c'_2, a'_1, a'_2)$ is a
specialization and $a'_1=a'_2$, we say that the specialized 
curves
$C'_1,C'_2$ are tangent at $a'_1$. In general this notion is
not
intrinsic to $C'_1,C'_2,a'_1$, but rather depends on the 
choice of
$C_1,C_2$. It does not always coincide with the usual notion
if
$X$ is an algebraic curve (and sometimes yields a more 
fruitful
notion, for our purposes).

\subheading{\num{5.} The field structure} From the given 
two-dimensional family
of curves, we may obtain a one-dimensional family of curves 
passing
through a single point $p_0$ in the plane. These curves are 
considered to
be multivalued functions from $X$ to $X$. As such they can 
be composed.
One would like to identify two curves that are tangent at 
$p_0$;
intuitively, an equivalence class of curves corresponds to a
slope;
one then wants to show that composition gives a well-defined
operation, 
which corresponds to multiplication on the slopes and at all
events
gives a group structure. In practice, tangency is not 
necessarily an
equivalence relation, and a number of other technical 
problems arise.
However, one can find curves $c_{ij}$ such that 
$c_{ij}^{-1}c_{i'j}$
and $c_{ij'}^{-1}c_{'ij'}$ are tangent for all $i<i'$ and 
$j<j'$. One then
shows that this gives an array as in (3), and hence gives an
Abelian group
structure.

Eugenia Rabinovich observed that in some cases the group 
obtained in this way
is in fact the additive rather than the multiplicative 
group.

One would like to find the field directly, by considering 
functions from the
plane to itself, modulo tangency, thus interpreting the 
tangent space to the
plane and what should be $\gl_2$ acting on it. This approach
poses severe
technical problems. Hence, one first interprets an Abelian 
group as above,
then uses it as a crutch to find the field configuration 
described above.
We find curves satisfying, up to tangency, the relations 
noted following
Proposition 4, in which addition is interpreted as the 
given Abelian
group structure and multiplication is interpreted as 
composition. We
show that this gives the field array.

\thm{\num{6.} Theorem 3} So far we have essentially 
described the
proof of Theorem \RM2. For Theorem \RM3 one is given a field
$F$ and
a Zariski geometry on $F$, refining the usual Zariski 
topology.
One must show that they coincide. If they do not, it can be 
shown that there
exists a plane curve $C$ not contained in any algebraic 
curve. From this
hypothesis one must obtain a contradiction.\ethm

We do this by developing an analog of Bezout's theorem. One 
defines the degree
of a curve $C$ to be the number of points of intersection of
 $C$ with a
generic line. Bezout's theorem then states that the number 
of points of
intersection of $C$ with an algebraic curve of degree $d$ is
at most
$d\boldcdot\deg(C)$. We give a version of a classical proof 
of this, 
``moving'' from an arbitrary algebraic curve to a generic 
one, and from there
to a special one for which the result is clear (the unions 
of $d$ lines).
A certain amount of preliminary work is required, showing 
that the projective
plane over $F$ is sufficiently complete for the purposes of 
such moves by
specializations.

Now one considers intersections of $C$ with algebraic curves
of high degree
$d$. The number of points of intersection increases linearly
with $d$, but the
dimension of the space of curves of degree $d$ increases 
quadratically; hence,
for large $d$ one can find a curve intersecting $C$ in more 
than
$d\boldcdot\deg(C)$ points. This is only possible if the 
intersection is
infinite and hence contains $C$, thus contradicting the 
choice of $C$ above.

Theorem 1 follows from Theorem 2, Theorem 3, and an 
analysis of covers in
the Zariski category. This analysis shows that any ample 
Zariski geometry
$X$ is a cover of a canonical algebraic curve $C$, such that
the pullback of a
curve on $C^n$ is typically irreducible in $X^n$. It follows
that if $X$ is
very ample, then $X=C$. It may be in general that the 
richness of the geometry
of $X$ arises entirely from that of $C$, but we do not know 
a precise
statement.

\heading Acknowledgment\endheading

The authors would like to thank Gregory Cherlin, Moshe 
Jarden, and Mike Fried
for close readings of this manuscript.

\enddocument